\documentclass[a4paper]{article}
\usepackage{lmodern}
\pdfoutput=1
\usepackage{graphicx,wrapfig}
\usepackage{lipsum}
\usepackage{amsmath,amssymb,mathrsfs}
\usepackage{hyperref}
\usepackage{url}
\usepackage{mathtools}
\usepackage{changepage}
\usepackage[ruled,vlined]{algorithm2e}
\usepackage{amsmath}
\usepackage{multirow}
\usepackage{wrapfig}
\usepackage{subfig}
\usepackage{color}
\usepackage{epsfig,enumerate,amsmath,amsfonts,amssymb,amsthm,mathrsfs,ifpdf}
\usepackage{indentfirst,relsize}
\usepackage{setspace,graphicx}
\usepackage{latexsym}
\usepackage[all]{xy}
\usepackage[dvipsnames]{pstricks}
\usepackage{pst-grad} 
\usepackage{pst-plot} 
\usepackage[margin = 2.50cm]{geometry}

\usepackage{authblk} 
\usepackage{algpseudocode}
\usepackage{color,soul}
\usepackage{enumitem}

\usepackage{orcidlink}
\usepackage{tabularx}
\usepackage{pst-grad} 
\usepackage{pst-plot} 
\usepackage{enumitem}

\usepackage{orcidlink}

\usepackage{tikz}
\usetikzlibrary{arrows.meta,positioning}

\newcommand{\remove}[1]{}

\newtheorem{theorem}{Theorem}
\newtheorem{lemma}[theorem]{Lemma}
\newtheorem{proposition}[theorem]{Proposition}
\newtheorem{corollary}[theorem]{Corollary}

\newtheorem{remark}{Remark}

\usetikzlibrary{arrows.meta,positioning,fit}
\usepackage[ruled,vlined]{algorithm2e}

\SetKwInput{KwInput}{Input}
\SetKwInput{KwOutput}{Output}

\SetKwInput{KwInput}{Input}
\SetKwInput{KwOutput}{Output}

\DeclareMathOperator{\Tr}{Tr}

\usepackage{todonotes}
\usepackage[most]{tcolorbox}

\title{On the Complexity and Algorithms for the Upper Domatic Number of Graphs}

\author[1]{Subhabrata Paul\footnote{subhabrata@iitp.ac.in}}

\author[2]{Kamal Santra \orcidlink{0009-0006-5997-1452} \footnote{kamal.7.2013@gmail.com, kamal.santra@iitg.ac.in}}

\affil[1]{Department of Mathematics\\
	
	Indian Institute of Technology Patna\\
	
	Patna, 801106, Bihar, India}

\affil[2]{Department of Mathematics\\
	
	Indian Institute of Technology Guwahati\\
	
	Guwahati, 781039, Assam, India}

\date{}

\begin{document}

	\maketitle
	\begin{abstract}
		Let \(G\) be a graph. For two disjoint vertex sets \(A,B\subseteq V(G)\), we say that \(A\) dominates \(B\) if every vertex of \(B\) has a neighbour in \(A\). An upper domatic partition of \(G\) is a partition \(\pi=\{V_1,V_2,\ldots,V_k\}\) of \(V(G)\) such that, for every two distinct parts \(V_i\) and \(V_j\), either \(V_i\) dominates \(V_j\), or \(V_j\) dominates \(V_i\), or both. The maximum order of such a partition is the upper domatic number of \(G\), denoted by \(D(G)\). In this paper, we study the computational complexity of the upper domatic number from both hardness and algorithmic perspectives. Motivated by the complexity questions raised in Phillips's thesis, we first prove that deciding whether \(D(G)\geq k\) is NP-complete when \(k\) is part of the input. Our reduction is from \textsc{Clique}, and it also shows NP-completeness for connected graphs of diameter two. On the positive side, we give exact algorithms for several graph classes. We prove that, for cographs, \(D(G)=\Tr(G)\), and we provide an \(O(n^3)\)-time cotree dynamic program. For unicyclic graphs, we use the known equality \(D(G)=\Tr(G)\). By deleting an edge of the unique cycle, we reduce the computation of \(\Tr(G)\) to the transitivity of a tree and one additional decision problem testing whether adding back the deleted edge increases the transitivity by one. This yields an \(O(n^3)\)-time algorithm for computing the transitivity, and hence the upper domatic number, of unicyclic graphs. We further prove that \(D(G)=\Tr(G)\) for complements of bipartite graphs, which yields a linear-time algorithm for complements of bipartite chain graphs. Finally, we show that \(D(G)=\Tr(G)\) for split graphs and obtain a linear-time algorithm for this class.
	\end{abstract}
	


	\section{Introduction}
\label{sec:introduction}

Graph partitioning is one of the central themes in graph theory. Many classical graph parameters are defined by partitioning the vertex set of a graph into parts satisfying some prescribed structural condition. For example, a proper coloring partitions the vertex set into independent sets, while a domatic partition partitions the vertex set into dominating sets. Such parameters are important not only because they measure structural features of graphs, but also because they often lead to natural algorithmic and complexity questions.

The domatic number was introduced by Cockayne and Hedetniemi~\cite{cockayne1977towards}. A domatic partition of a graph \(G\) is a partition of \(V(G)\) into dominating sets, and the maximum order of such a partition is the domatic number \(d(G)\). Around the same line of study, Christen and Selkow introduced Grundy colorings and the Grundy number~\cite{christen1979some}. A Grundy partition is a proper coloring with an additional ordering condition: every vertex in a later color class must have a neighbour in each earlier color class. Harary and Hedetniemi studied the achromatic number, which is another extremal coloring parameter based on requiring an edge between every pair of color classes~\cite{harary1970achromatic}. These parameters show how graph coloring and domination-type conditions naturally interact through vertex partitions.

A useful way to view these partition parameters is through domination between parts. Let \(G\) be a graph, and let \(A,B\subseteq V(G)\) be two disjoint vertex sets. We say that \(A\) dominates \(B\), written \(A\to B\), if every vertex of \(B\) has a neighbour in \(A\). This type of domination relation between parts is closely connected with domination digraphs~\cite{goddard2012algorithmic}. Hedetniemi and Hedetniemi introduced the transitivity of a graph through ordered domination partitions~\cite{hedetniemi2018transitivity}. More precisely, a transitive partition is an ordered partition \((V_1,V_2,\ldots,V_k)\) of \(V(G)\) such that \(V_i\to V_j\) whenever \(i<j\). The maximum order of such a partition is the transitivity of \(G\), denoted by \(\operatorname{Tr}(G)\).

The study of transitivity is also related to earlier work on iterated colorings of graphs, where domination-type and independence-type partition parameters were considered from an algorithmic point of view~\cite{hedetniemi2004iterated}. Further results on transitivity were obtained for special graph classes, including cactus graphs, coronas, Cartesian products, and joins~\cite{haynes2019transitivity}. More recently, Paul and Santra studied the maximum transitivity problem on subclasses of bipartite graphs, proving hardness for perfect elimination bipartite graphs and giving linear-time algorithms for bipartite chain graphs~\cite{paulSantra2023transitivityBipartite}. They also studied transitivity for split graphs, pseudo-split graphs, complements of bipartite chain graphs, central graphs, and transitively critical graphs~\cite{paulSantra2025transitivityGraphClasses}. Since Grundy partitions require the same ordered domination condition together with the additional requirement that each part is independent~\cite{christen1979some}, transitive partitions may be viewed as a relaxation of Grundy partitions.

The upper domatic number was introduced by Haynes et al.~\cite{haynes2020upper}. An upper domatic partition of \(G\) is a partition \(\pi=\{V_1,V_2,\ldots,V_k\}\) of \(V(G)\) such that, for every two distinct parts \(V_i\) and \(V_j\), either \(V_i\to V_j\), or \(V_j\to V_i\), or both. The maximum order of an upper domatic partition is the upper domatic number of \(G\), denoted by \(D(G)\). An upper domatic partition of size \(D(G)\) is called a \(D\)-partition. In this sense, the upper domatic number is a relaxation of the domatic number: in a domatic partition, every pair of parts dominates each other, while in an upper domatic partition, domination is required in at least one direction.

Domination digraphs provide a convenient language for these notions. Given a vertex partition \(\pi=\{V_1,V_2,\ldots,V_k\}\) of \(V(G)\), the domination digraph \(D_G(\pi)\) has vertex set \(\pi\), and it contains an arc \(V_i\to V_j\) whenever \(V_i\) dominates \(V_j\) in \(G\). Domination digraphs were studied by Goddard et al.~\cite{goddard2012algorithmic}. In this language, a domatic partition corresponds to a bi-directed complete domination digraph, a transitive partition corresponds to a domination digraph containing a spanning transitive tournament, and an upper domatic partition corresponds to a complete domination digraph.

These definitions give the basic inequalities \(d(G)\leq \Tr(G)\leq D(G)\) for every graph \(G\). Moreover, every upper domatic partition gives a pseudoachromatic coloring, because every two parts have at least one edge between them; hence \(D(G)\leq \psi(G)\), where \(\psi(G)\) denotes the pseudoachromatic number. Thus, the upper domatic number lies naturally between domination-type and complete-coloring-type partition parameters. This makes it a useful parameter for understanding how domination relations between color classes influence the difficulty of graph partitioning problems.

Haynes et al.~\cite{haynes2020upper} established several fundamental properties of \(D(G)\). They proved the upper bound \(D(G)\leq \Delta(G)+1\), where \(\Delta(G)\) is the maximum degree of \(G\), and studied the relationship between \(D(G)\) and \(\Tr(G)\). In particular, they showed that the two parameters coincide for several basic cases, including trees. Samuel and Joseph further developed the theory and proved that \(D(G)=\Tr(G)\) for unicyclic graphs~\cite{samuel2020new}. These results suggest that transitivity can often be used as a bridge for computing the upper domatic number. However, this equality does not hold for all graphs; therefore, it is important to identify graph classes where the equality is valid and algorithmically useful.

The computational complexity of the upper domatic number was already raised as a natural question in the early study of the parameter. In his master's thesis, Phillips listed as an open problem the time complexity of deciding whether a graph satisfies \(D(G)\geq 4\) or \(D(G)\geq 5\)~\cite{phillips2017study}. This motivates the complexity part of the present paper. We prove NP-completeness of the general decision version, where the threshold \(k\) is part of the input. This result does not by itself settle the fixed-threshold questions \(D(G)\geq 4\) and \(D(G)\geq 5\), but it gives a complete hardness result for the natural variable-threshold formulation.

In this paper, we study the upper domatic number from two complementary viewpoints. First, we prove that the general decision problem is computationally hard. Second, we give polynomial-time algorithms for several structured graph classes. A common theme in the positive results is to prove \(D(G)=\Tr(G)\), and then compute the transitivity using structural properties of the graph class.

Our contributions are as follows.

\begin{enumerate}
	\item We prove that the decision problem asking whether \(D(G)\geq k\) is NP-complete. The reduction is from the classical \textsc{Clique} problem~\cite{garey1979computers}. Given an instance \((G,q)\), we construct a graph \(H\) by taking a copy of \(G\), adding an independent set of size \(n-q+1\), and joining every new vertex to every old vertex. A clique of size \(q\) in \(G\) gives an upper domatic partition of \(H\) with \(n+1\) parts. Conversely, any upper domatic partition of \(H\) with at least \(n+1\) parts must contain many singleton parts, and these singleton parts force a clique in the original graph. We also show that the construction can be modified so that the resulting graph is connected and has diameter exactly two.
	
	\item We study cographs. We prove that \(D(G)=\Tr(G)\) for every cograph \(G\). This is not immediate, because \(D(G)\) does not satisfy a simple equality rule over disconnected components. We use the disconnected upper-domatic bound of Samuel and Joseph~\cite{samuel2020new}, together with the component formula for transitivity, to handle union nodes. For join nodes, we show that mixed parts can be placed at the beginning of a transitive partition. After proving the equality \(D(G)=\Tr(G)\), we compute \(\Tr(G)\) by a dynamic program over a binary cotree. Using the cotree structure of cographs~\cite{corneil1981complement,corneil1985linear}, this gives an \(O(n^3)\)-time algorithm.
	
	\item We give an \(O(n^3)\)-time algorithm for unicyclic graphs. Since \(D(G)=\Tr(G)\) in this class, it is enough to compute the transitivity. We delete an arbitrary edge \(xy\) of the unique cycle, obtaining a tree \(T=G-xy\). By the edge-addition and edge-deletion bounds of Hedetniemi and Hedetniemi~\cite{hedetniemi2018transitivity}, we have \(\Tr(G)\in\{\Tr(T),\Tr(T)+1\}\). Since \(\Tr(T)\) can be computed in linear time using the known tree algorithm, it remains only to decide whether adding back the edge \(xy\) increases the transitivity by one. We solve this decision problem by a dynamic program on the rooted tree \(T\), using a transitive-labelling formulation and explicitly accounting for the single restored edge \(xy\). This yields an \(O(n^3)\)-time algorithm for computing both \(\Tr(G)\) and \(D(G)\) on unicyclic graphs.
	
	\item We prove that \(D(G)=\Tr(G)\) for complements of bipartite graphs. The key step is to show that every such graph admits a maximum upper domatic partition containing a source set. Once such a source set is found, it can be placed first, and the argument proceeds by induction on the remaining graph. As a consequence, using the known linear-time algorithm for transitivity on complements of bipartite chain graphs~\cite{paulSantra2025transitivityGraphClasses}, we obtain a linear-time algorithm for the upper domatic number on complements of bipartite chain graphs.
	
	\item We study split graphs. If \(G\) is a split graph with clique number \(\omega(G)\), then we prove that \(D(G)\) is either \(\omega(G)\) or \(\omega(G)+1\). More precisely, the value depends on whether every vertex of a maximum clique has a neighbour in the independent set. We then prove that \(D(G)=\Tr(G)\) for split graphs. Together with the known linear-time algorithm for transitivity on split graphs~\cite{paulSantra2025transitivityGraphClasses}, this yields a linear-time algorithm for computing the upper domatic number on split graphs.
\end{enumerate}

The paper is organized as follows. Section~\ref{sec:preliminaries} gives the notation and basic definitions. Section~\ref{sec:upper-domatic-np-complete} proves NP-completeness of the decision problem. Section~\ref{sec:upper-domatic-cographs} studies cographs and gives the cotree dynamic program. Section~\ref{sec:unicyclic} gives the algorithm for unicyclic graphs. Section~\ref{sec:complement-bipartite} proves the equality \(D(G)=\Tr(G)\) for complements of bipartite graphs and derives the algorithmic consequence for complements of bipartite chain graphs. Section~\ref{sec:split-graphs} treats split graphs. Finally, Section~\ref{sec:conclusion} concludes the paper.

\section{Notation and Definitions}
\label{sec:preliminaries}

Throughout this paper, all graphs are finite, simple, and undirected. Let \(G=(V,E)\) be a graph. The order of \(G\) is \(|V(G)|\), and the size of \(G\) is \(|E(G)|\). When the graph is clear from the context, we write \(V\) and \(E\) instead of \(V(G)\) and \(E(G)\), respectively. For two vertices \(u,v\in V(G)\), we write \(uv\in E(G)\) to mean that \(u\) and \(v\) are adjacent.

For a vertex \(v\in V(G)\), the open neighbourhood of \(v\) is denoted by \(N_G(v)\), and the closed neighbourhood of \(v\) is \(N_G[v]=N_G(v)\cup\{v\}\). The degree of \(v\) is \(\deg_G(v)=|N_G(v)|\). The maximum and minimum degrees of \(G\) are denoted by \(\Delta(G)\) and \(\delta(G)\), respectively. When there is no chance of confusion, we simply write \(N(v)\), \(N[v]\), and \(\deg(v)\). For a set \(S\subseteq V(G)\), we write \(N_G(S)=\bigcup_{v\in S}N_G(v)\) and \(N_G[S]=N_G(S)\cup S\).

A graph \(H=(V',E')\) is a subgraph of \(G\) if \(V'\subseteq V(G)\) and \(E'\subseteq E(G)\). If \(S\subseteq V(G)\), then the subgraph induced by \(S\) is denoted by \(G[S]\). Thus, \(G[S]\) has vertex set \(S\) and contains exactly those edges of \(G\) whose two endpoints both belong to \(S\). We write \(G-S\) for the induced subgraph \(G[V(G)\setminus S]\). For a vertex \(v\), we write \(G-v\) instead of \(G-\{v\}\). The complement of \(G\) is denoted by \(\overline{G}\); it has the same vertex set as \(G\), and two distinct vertices are adjacent in \(\overline{G}\) exactly when they are nonadjacent in \(G\).

A set \(S\subseteq V(G)\) is independent if no two vertices of \(S\) are adjacent. A set \(K\subseteq V(G)\) is a clique if every two distinct vertices of \(K\) are adjacent. The maximum size of a clique in \(G\) is the clique number of \(G\), denoted by \(\omega(G)\). A vertex of degree one is called a leaf. A graph is connected if every two vertices are joined by a path. A connected graph containing exactly one cycle is called unicyclic. Equivalently, a connected graph on \(n\) vertices is unicyclic if it has exactly \(n\) edges.

We shall use two basic graph operations. The disjoint union of two vertex-disjoint graphs \(G_1\) and \(G_2\), denoted by \(G_1\cup G_2\), is the graph with vertex set \(V(G_1)\cup V(G_2)\) and edge set \(E(G_1)\cup E(G_2)\). The join of \(G_1\) and \(G_2\), denoted by \(G_1\vee G_2\), is obtained from \(G_1\cup G_2\) by adding all possible edges between \(V(G_1)\) and \(V(G_2)\).

A set \(D\subseteq V(G)\) is a dominating set of \(G\) if every vertex in \(V(G)\setminus D\) has a neighbour in \(D\). For two disjoint sets \(A,B\subseteq V(G)\), we say that \(A\) dominates \(B\), written \(A\to B\), if every vertex of \(B\) has a neighbour in \(A\). Notice that this is an open-neighbourhood condition between two disjoint sets: vertices of \(B\) must be dominated by vertices of \(A\), and vertices do not dominate themselves in this relation.

A partition of \(V(G)\) is a family \(\pi=\{V_1,V_2,\ldots,V_k\}\) of nonempty, pairwise disjoint subsets of \(V(G)\) whose union is \(V(G)\). The sets \(V_1,V_2,\ldots,V_k\) are called the parts of the partition. A partition is ordered if an order of its parts is fixed, and in that case we write it as \((V_1,V_2,\ldots,V_k)\).

A domatic partition of \(G\) is a partition of \(V(G)\) into dominating sets. The maximum order of a domatic partition is the domatic number of \(G\), denoted by \(d(G)\), introduced by Cockayne and Hedetniemi~\cite{cockayne1977towards}. A Grundy partition is an ordered partition \((V_1,V_2,\ldots,V_k)\) such that each \(V_i\) is independent and, for every \(1\leq i<j\leq k\), the part \(V_i\) dominates \(V_j\). The maximum order of a Grundy partition is the Grundy number of \(G\), denoted by \(\Gamma(G)\)~\cite{christen1979some}.

A transitive partition of \(G\) is an ordered partition \((V_1,V_2,\ldots,V_k)\) such that \(V_i\to V_j\) for every \(1\leq i<j\leq k\). The maximum order of a transitive partition is the transitivity of \(G\), denoted by \(\Tr(G)\)~\cite{hedetniemi2018transitivity}. Thus, a transitive partition is obtained from a Grundy partition by dropping the requirement that the parts must be independent.

An upper domatic partition of \(G\) is a partition \(\pi=\{V_1,V_2,\ldots,V_k\}\) such that, for every two distinct parts \(V_i\) and \(V_j\), either \(V_i\to V_j\), or \(V_j\to V_i\), or both. The maximum order of an upper domatic partition is the upper domatic number of \(G\), denoted by \(D(G)\)~\cite{haynes2020upper}. An upper domatic partition of size \(D(G)\) is called a \(D\)-partition of \(G\).

The above definitions immediately imply the inequalities \(\Gamma(G)\leq \Tr(G)\leq D(G)\). Also, every domatic partition is an upper domatic partition, and hence \(d(G)\leq D(G)\). In fact, every domatic partition is also a transitive partition in any ordering of its parts, so \(d(G)\leq \Tr(G)\leq D(G)\).

It is often useful to encode the domination relations between parts by a directed graph. Given a partition \(\pi=\{V_1,V_2,\ldots,V_k\}\) of \(V(G)\), the domination digraph of \(\pi\), denoted by \(D_G(\pi)\), is the digraph with vertex set \(\pi\), where there is an arc \(V_i\rightarrow V_j\) if and only if \(V_i\to V_j\) in \(G\). Domination digraphs were studied by Goddard et al.~\cite{goddard2012algorithmic}. A digraph is complete if, for every two distinct vertices \(x\) and \(y\), at least one of the arcs \(x\to y\) and \(y\to x\) is present. It is bi-directed complete if, for every two distinct vertices \(x\) and \(y\), both arcs \(x\to y\) and \(y\to x\) are present. A tournament is a complete digraph in which exactly one arc is present between every pair of distinct vertices. A tournament is transitive if its vertices can be ordered as \(x_1,x_2,\ldots,x_k\) such that \(x_i\to x_j\) whenever \(i<j\).

With this terminology, a partition \(\pi\) is an upper domatic partition exactly when \(D_G(\pi)\) is complete. It is a domatic partition exactly when \(D_G(\pi)\) is bi-directed complete. An ordered partition \((V_1,V_2,\ldots,V_k)\) is a transitive partition exactly when \(D_G(\pi)\) contains all arcs \(V_i\to V_j\) with \(i<j\); equivalently, \(D_G(\pi)\) contains a spanning transitive tournament.

Let \(\pi=\{V_1,V_2,\ldots,V_k\}\) be a partition of \(V(G)\). A part \(V_i\) is called a source set of \(\pi\) if \(V_i\) dominates every other part of \(\pi\). Equivalently, \(V_i\) has out-degree \(k-1\) in the domination digraph \(D_G(\pi)\). Similarly, \(V_i\) is called a sink set of \(\pi\) if every other part of \(\pi\) dominates \(V_i\). Equivalently, \(V_i\) has in-degree \(k-1\) in \(D_G(\pi)\). In every transitive partition, the first part is a source set and the last part is a sink set.

We also recall the graph classes considered in this paper. A cograph is a graph obtained from single vertices by repeated applications of disjoint union and join. Equivalently, cographs are precisely the \(P_4\)-free graphs, and every cograph admits a cotree representation whose leaves correspond to the vertices and whose internal nodes correspond to union and join operations~\cite{corneil1981complement}. A cotree can be constructed in linear time~\cite{corneil1985linear}.

A split graph is a graph whose vertex set can be partitioned into a clique \(K\) and an independent set \(S\). Such a partition is called a split partition. In the split-graph section, we shall use a split partition \(V(G)=K\cup S\) where \(K\) is a maximum clique, so \(|K|=\omega(G)\).

A graph is bipartite if its vertex set can be partitioned into two independent sets \(X\) and \(Y\). A graph \(G\) is the complement of a bipartite graph if \(V(G)\) can be partitioned into two cliques \(X\) and \(Y\). Indeed, if \(\overline{G}\) is bipartite with bipartition \(X\cup Y\), then both \(X\) and \(Y\) are cliques in \(G\). A bipartite graph is a bipartite chain graph if the vertices in each part can be ordered so that their neighbourhoods are nested. Equivalently, bipartite chain graphs are the \(2K_2\)-free bipartite graphs~\cite{yannakakis1982complexity}. Linear-time recognition and construction of such orderings are known~\cite{heggernes2007linear}.

Finally, a unicyclic graph is a connected graph with exactly one cycle. In the unicyclic-graph section, we repeatedly use the fact that deleting any edge of the unique cycle produces a tree. This allows us to root the resulting tree and process vertices in a postorder.

\section{NP-completeness of the Upper Domatic Number Problem}
\label{sec:upper-domatic-np-complete}

In this section, we prove that deciding whether the upper domatic number is at least a given integer is NP-complete. The complexity of threshold versions of this problem was already raised by Phillips, who asked about the time complexity of deciding whether \(D(G)\geq 4\) or \(D(G)\geq 5\)~\cite{phillips2017study}. Here, we prove the NP-completeness of the general decision problem, where the threshold is part of the input. We use a reduction from the classical \textsc{Clique} problem, which is NP-complete~\cite{garey1979computers}. The main idea is simple. Starting from an instance \((G,q)\) of \textsc{Clique}, we build a graph \(H\) by taking a copy of \(G\), adding an independent set of carefully chosen size, and joining every new vertex to every old vertex. If \(G\) has a clique of size \(q\), then \(H\) has a large upper domatic partition. Conversely, any sufficiently large upper domatic partition of \(H\) must contain many singleton parts, and these singleton parts force a clique in \(G\).

\begin{figure}[htbp]
	\centering
	\resizebox{0.75\linewidth}{!}{%
		\begin{tikzpicture}[
			vertex/.style={circle,draw,fill=white,inner sep=1.2pt,minimum size=5.5mm},
			yvertex/.style={circle,draw,fill=gray!20,inner sep=1.2pt,minimum size=5.5mm},
			edge/.style={line width=0.45pt},
			joinedge/.style={line width=0.35pt,dashed,gray},
			box/.style={draw,rounded corners,inner sep=5pt},
			lab/.style={font=\small},
			>=Latex
			]
			
			\node[vertex] (g1) at (0,1.0) {};
			\node[vertex] (g2) at (1.0,1.2) {};
			\node[vertex] (g3) at (1.4,0.2) {};
			\node[vertex] (g4) at (0.4,-0.55) {};
			\draw[edge] (g1)--(g2);
			\draw[edge] (g2)--(g3);
			\draw[edge] (g3)--(g4);
			\draw[edge] (g4)--(g1);
			\draw[edge] (g1)--(g3);
			\node[box,fit=(g1)(g2)(g3)(g4)] (Gbox) {};
			\node[lab,below=4pt of Gbox] {Input graph \(G\)};
			
			\node at (2.75,0.25) {\(\Longrightarrow\)};
			
			\node[vertex] (x1) at (4.2,1.0) {};
			\node[vertex] (x2) at (5.2,1.2) {};
			\node[vertex] (x3) at (5.6,0.2) {};
			\node[vertex] (x4) at (4.6,-0.55) {};
			\draw[edge] (x1)--(x2);
			\draw[edge] (x2)--(x3);
			\draw[edge] (x3)--(x4);
			\draw[edge] (x4)--(x1);
			\draw[edge] (x1)--(x3);
			\node[box,fit=(x1)(x2)(x3)(x4)] (Xbox) {};
			\node[lab,below=4pt of Xbox] {Copy \(X\cong G\)};
			
			\node[yvertex] (y1) at (8.7,1.0) {};
			\node[yvertex] (y2) at (8.7,0.25) {};
			\node[yvertex] (y3) at (8.7,-0.50) {};
			\node[box,fit=(y1)(y2)(y3)] (Ybox) {};
			\node[lab,below=4pt of Ybox] {Independent set \(Y\)};
			
			\foreach \x in {x1,x2,x3,x4}{
				\foreach \y in {y1,y2,y3}{
					\draw[joinedge] (\x)--(\y);
				}
			}
			
			\node[lab,above=8pt of Xbox, xshift=1.9cm] {Construct \(H\)};
			\node[lab,above=1pt of Ybox] {\(|Y|=n-q+1\)};
			
		\end{tikzpicture}%
	}
	\caption{The reduction from \textsc{Clique} to the upper domatic number problem. The graph \(H\) is obtained by taking a copy \(X\cong G\), adding an independent set \(Y\) of size \(n-q+1\), and adding all edges between \(X\) and \(Y\). The dashed edges indicate this complete join.}
	\label{fig:clique-to-upper-domatic}
\end{figure}

We consider the following decision problem.

\smallskip
\noindent\textsc{Upper Domatic Number}

\noindent\emph{Instance:} A graph \(G\) and an integer \(k\).

\noindent\emph{Question:} Is \(D(G)\geq k\)?

\begin{theorem}
	\label{thm:upper-domatic-np-complete}
	The \textsc{Upper Domatic Number} decision problem is NP-complete.
\end{theorem}

\begin{proof}
	The problem belongs to NP. A certificate is a partition \(\mathcal P=\{P_1,\ldots,P_r\}\) of \(V(G)\) into \(r\geq k\) nonempty parts. For each pair of parts \(P_i\) and \(P_j\), we can check in polynomial time whether \(P_i\) dominates \(P_j\) or \(P_j\) dominates \(P_i\). Hence, the certificate can be verified in polynomial time.
	
	We prove NP-hardness by reducing from \textsc{Clique}. Let \((G,q)\) be an instance of \textsc{Clique}, where \(G=(V,E)\), \(|V|=n\), and the question is whether \(G\) contains a clique of size at least \(q\). We may assume \(1\leq q\leq n\). We construct a graph \(H\) as follows. First, take a copy \(X=\{x_v:v\in V\}\) of \(V(G)\), and put \(x_u x_v\in E(H)\) exactly when \(uv\in E(G)\). Thus, \(H[X]\) is isomorphic to \(G\). Next, add an independent set \(Y\) of size \(n-q+1\). Finally, join every vertex of \(Y\) to every vertex of \(X\). There are no edges inside \(Y\). The construction is shown schematically in Figure~\ref{fig:clique-to-upper-domatic}. We set the target value to be \(K=n+1\). The graph \(H\) has \(2n-q+1\leq 2n\) vertices and \(|E(G)|+n(n-q+1)\leq |E(G)|+n^2\) edges, so the construction is polynomial.
	
	We claim that \(G\) has a clique of size \(q\) if and only if \(D(H)\geq n+1\). First, suppose that \(G\) has a clique \(C\) of size \(q\). Let \(X_C=\{x_v:v\in C\}\). Write \(X\setminus X_C=\{x_1,\ldots,x_{n-q}\}\), and write \(Y=\{y_0,y_1,\ldots,y_{n-q}\}\). We define a partition of \(V(H)\) as follows. For each vertex \(x\in X_C\), make the singleton part \(\{x\}\). For each \(i\in\{1,\ldots,n-q\}\), make the paired part \(\{x_i,y_i\}\). Finally, make the singleton part \(\{y_0\}\). This gives \(q+(n-q)+1=n+1\) parts.
	
	We now verify that this is an upper domatic partition. The singleton vertices in \(X_C\) are pairwise adjacent because \(C\) is a clique in \(G\). The singleton \(\{y_0\}\) is comparable with every singleton \(\{x\}\), where \(x\in X_C\), because \(y_0\) is adjacent to every vertex of \(X\). A paired part \(\{x_i,y_i\}\) dominates a singleton \(\{x\}\) with \(x\in X_C\), since \(y_i\) is adjacent to \(x\). Also, \(\{x_i,y_i\}\) dominates \(\{y_0\}\), since \(x_i\) is adjacent to \(y_0\). Finally, if \(i\neq j\), then \(\{x_i,y_i\}\) dominates \(\{x_j,y_j\}\), because \(y_i\) is adjacent to \(x_j\), and \(x_i\) is adjacent to \(y_j\). Thus, every pair of parts is comparable by domination, and so \(D(H)\geq n+1\).
	
	Conversely, suppose that \(D(H)\geq n+1\). Then \(H\) has an upper domatic partition with \(r\geq n+1\) parts. The graph \(H\) has \(|V(H)|=n+(n-q+1)=2n-q+1\) vertices. In any partition of an \(N\)-vertex set into \(r\) nonempty parts, the number of singleton parts is at least \(2r-N\), because if \(s\) parts are singletons, then \(N\geq s+2(r-s)=2r-s\). Therefore, this upper domatic partition has at least \(2(n+1)-(2n-q+1)=q+1\) singleton parts.
	
	In any upper domatic partition, the vertices appearing in singleton parts induce a clique~\cite{haynes2020upper}. Indeed, if \(\{a\}\) and \(\{b\}\) are two singleton parts, then one of these singleton parts must dominate the other, and this is possible only when \(a\) and \(b\) are adjacent. Hence, \(H\) contains a clique of size at least \(q+1\) formed by singleton vertices. Since \(Y\) is independent, such a clique contains at most one vertex from \(Y\). Therefore, at least \(q\) of these singleton vertices lie in \(X\). But \(H[X]\) is isomorphic to \(G\), and so \(G\) contains a clique of size at least \(q\).
	
	Thus, \((G,q)\) is a yes-instance of \textsc{Clique} if and only if \((H,n+1)\) is a yes-instance of \textsc{Upper Domatic Number}. Therefore, the decision version of the upper domatic number problem is NP-hard. Since the problem is also in NP, it is NP-complete.
\end{proof}

\subsection{The graph class obtained in the reduction}
\label{subsec:reduced-class-upper-domatic}

The graph produced by the reduction has the form \(H=G\vee\overline{K}_{n-q+1}\). Equivalently, its vertex set admits a partition \(V(H)=X\cup Y\) such that \(H[X]\cong G\), the set \(Y\) is independent, and every vertex of \(Y\) is adjacent to every vertex of \(X\). In addition, \(1\leq |Y|\leq |X|\). All vertices of \(Y\) have the same open neighbourhood \(X\), so \(Y\) is an independent module consisting of false twins.

The construction does not restrict the output to split graphs, chordal graphs, block graphs, or cographs. Indeed, \(H[X]\cong G\) is an induced subgraph of \(H\), while \(G\) is arbitrary. For example, if \(G=C_4\), then \(H\) contains an induced \(C_4\), and hence, \(H\) need not be chordal, split, or a block graph. Similarly, taking \(G=P_4\) shows that the resulting graph need not be a cograph.

Nevertheless, every graph produced by the reduction is connected and has diameter at most two. Every vertex of \(X\) is adjacent to every vertex of \(Y\); two nonadjacent vertices of \(X\) have a common neighbour in \(Y\); and two vertices of \(Y\) have a common neighbour in \(X\). To obtain a diameter exactly two in every instance, we apply a small preprocessing step to the \textsc{Clique} instance.

\begin{corollary}
	\label{cor:upper-domatic-diameter-two}
	The \textsc{Upper Domatic Number} decision problem is NP-complete even when restricted to connected graphs of diameter exactly two.
\end{corollary}

\begin{proof}
	Membership in NP follows from Theorem~\ref{thm:upper-domatic-np-complete}. For NP-hardness, let \((G,q)\) be an instance of \textsc{Clique}, where \(|V(G)|=n\). We may assume \(q\geq 2\), since the cases \(q\leq 1\) are trivial. Add one isolated vertex \(z\) to \(G\), and denote the resulting graph by \(G^{+}\). Let \(n^{+}=|V(G^{+})|=n+1\). Since \(q\geq 2\), the isolated vertex \(z\) cannot belong to a clique of size \(q\), and hence \(G\) has a clique of size at least \(q\) if and only if \(G^{+}\) does.
	
	Apply the reduction of Theorem~\ref{thm:upper-domatic-np-complete} to the instance \((G^{+},q)\). Thus, construct \(H=G^{+}\vee\overline{K}_{n^{+}-q+1}\) and set the target value to \(n^{+}+1\). By the correctness of the reduction, \(G\) has a clique of size at least \(q\) if and only if \(D(H)\geq n^{+}+1\).
	
	The graph \(H\) is connected, and every two vertices are at distance at most two. Moreover, \(z\) is nonadjacent in \(G^{+}\) to every original vertex of \(G\). If \(x\in V(G)\) and \(y\) is any vertex of the added independent set, then \(z-y-x\) is a path of length two in \(H\). Therefore, \(d_H(z,x)=2\), and consequently the diameter of \(H\) is exactly two. Hence, the problem remains NP-hard on connected graphs of diameter two.
\end{proof}

\begin{remark}
	The reduction shows NP-completeness for the class of graphs admitting a partition \(V(H)=X\cup Y\), where \(Y\) is an independent false-twin module complete to \(X\), and \(H[X]\) is arbitrary. Since this class is contained in the class of connected graphs of diameter at most two, and the preprocessing above ensures diameter exactly two, we obtain the more standard statement that the problem is NP-complete for connected graphs of diameter two.
\end{remark}

\section{Upper Domatic Partitions in Cographs}
\label{sec:upper-domatic-cographs}

In this section, we prove that the upper domatic number of a cograph can be computed in polynomial time. The main point is that we do not compute the upper domatic number directly at union nodes. This distinction is important because the upper domatic number does not behave by equality over disconnected components. Instead, we first prove that \(D(G)=\Tr(G)\) for cographs, and then compute \(\Tr(G)\) by dynamic programming over a cotree.

We shall use the following two facts about disconnected graphs. If \(G\) has connected components \(G_1,\ldots,G_r\), then Samuel and Joseph proved that \(D(G)\leq \max\{D(G_i):1\leq i\leq r\}\), and they gave an example showing that this inequality can be strict~\cite{samuel2020new}. On the other hand, transitivity satisfies the exact formula \(\Tr(G)=\max\{\Tr(G_i):1\leq i\leq r\}\) for disconnected graphs~\cite{hedetniemi2018transitivity,samuel2020new}. This is the reason why the dynamic program below is stated for \(\Tr\), not directly for \(D\).

Recall that cographs are precisely the graphs obtained from single vertices by repeated disjoint union and join operations. Equivalently, every cograph has a cotree whose leaves are the vertices and whose internal nodes are labelled by union or join~\cite{corneil1981complement}. A cotree can be constructed in linear time~\cite{corneil1985linear}. We use a binary cotree. If a cotree node has more than two children, we replace it with a binary tree of nodes with the same label; this does not change the represented graph and increases the number of nodes only linearly.

The following lemma is the key structural ingredient of this section. Although every transitive partition is an upper domatic partition, the converse is not true for arbitrary graphs. We show that, for cographs, the recursive structure given by union and join operations is strong enough to transform a maximum upper domatic partition into a transitive partition of the same size. Hence, on cographs, computing the upper domatic number is equivalent to computing transitivity.

\begin{lemma}
	\label{lem:cograph-D-equals-Tr}
	For every cograph \(G\), we have \(D(G)=\Tr(G)\).
\end{lemma}

\begin{proof}
	We prove the statement by induction on \(|V(G)|\). The result is immediate when \(|V(G)|=1\).
	
	First, assume that \(G\) is disconnected, and let \(C_1,\ldots,C_r\) be its connected components. Since every component of a cograph is again a cograph, the induction hypothesis gives \(D(C_i)=\Tr(C_i)\) for every \(i\). By the disconnected upper-domatic bound of Samuel and Joseph~\cite{samuel2020new}, we have
	\(
	D(G)\leq \max\{D(C_i):1\leq i\leq r\}.
	\)
	Therefore,
	\(
	D(G)\leq \max\{\Tr(C_i):1\leq i\leq r\}.
	\)
	By the component formula for transitivity~\cite{hedetniemi2018transitivity,samuel2020new}, the right-hand side is \(\Tr(G)\). Thus \(D(G)\leq \Tr(G)\). Since \(\Tr(G)\leq D(G)\) holds for every graph, it follows that \(D(G)=\Tr(G)\).
	
	Now assume that \(G\) is connected and has at least two vertices. Since \(G\) is a connected cograph, there exist smaller cographs \(G_1\) and \(G_2\) such that \(G=G_1\vee G_2\). Let \(\mathcal P\) be a \(D\)-partition of \(G\), and let \(k=|\mathcal P|=D(G)\). A part of \(\mathcal P\) is called mixed if it contains at least one vertex from \(G_1\) and at least one vertex from \(G_2\). Every mixed part dominates every other part: if the target vertex lies in \(G_1\), then the vertex of the mixed part lying in \(G_2\) is adjacent to it, and if the target vertex lies in \(G_2\), then the vertex of the mixed part lying in \(G_1\) is adjacent to it.
	
	Let \(X_i\) be the set of vertices of \(G_i\) that belong to mixed parts, for \(i\in\{1,2\}\). Consider the pure parts of \(\mathcal P\) contained in \(G_i\). If \(G_i-X_i\) is nonempty, then these pure parts form an upper domatic partition of \(G_i-X_i\). Indeed, domination between two pure parts on the same side cannot be helped by vertices from the other side, because both the source part and the target part are contained in the same graph \(G_i\). By the induction hypothesis, \(D(G_i-X_i)=\Tr(G_i-X_i)\). Hence, the pure parts on the \(G_i\)-side can be replaced by a transitive partition of \(G_i-X_i\) having at least the same number of parts. If \(G_i-X_i\) is empty, then there are no pure parts on that side and nothing needs to be replaced.
	
	We now build a transitive partition of \(G\) with at least \(k\) parts. Place all mixed parts first, in any order. Then append a transitive partition of \(G_1-X_1\), if \(G_1-X_1\) is nonempty, followed by a transitive partition of \(G_2-X_2\), if \(G_2-X_2\) is nonempty. The mixed parts dominate all later parts. Inside each nonempty side, transitivity holds by construction. Finally, every pure part on the \(G_1\)-side dominates every pure part on the \(G_2\)-side, because \(G\) is the join of \(G_1\) and \(G_2\). Thus, \(G\) has a transitive partition with at least \(k=D(G)\) parts. Therefore, \(\Tr(G)\geq D(G)\), and since \(\Tr(G)\leq D(G)\) holds for every graph, we obtain \(D(G)=\Tr(G)\).
\end{proof}

\subsection{Dynamic programming over the cotree}
\label{subsec:cograph-dp-recurrence}

Let \(T\) be a binary cotree of \(G\). For a node \(x\) of \(T\), let \(G_x\) be the cograph represented by the subtree rooted at \(x\), and let \(n_x=|V(G_x)|\). We use the convention that \(\Tr(\emptyset)=0\). For each integer \(s\), where \(0\leq s\leq n_x\), define
\(
F_x(s)=
\max\{\Tr(G_x-S):S\subseteq V(G_x),\ |S|=s\}.
\)
Thus, \(F_x(s)\) is the best transitivity value that can be obtained from \(G_x\) after deleting exactly \(s\) vertices. By Lemma~\ref{lem:cograph-D-equals-Tr}, the same table also gives the corresponding maximum upper domatic value among induced subgraphs obtained from \(G_x\) by deleting exactly \(s\) vertices.

If \(x\) is a leaf, then \(G_x\) consists of one vertex. Hence \(F_x(0)=1\), and after deleting that vertex, we obtain the empty graph, so \(F_x(1)=0\).

Now suppose that \(x\) is an internal node with children \(y\) and \(z\). If \(x\) is a union node, then \(G_x=G_y\cup G_z\). After deleting \(s_y\) vertices from \(G_y\) and \(s_z\) vertices from \(G_z\), where \(s_y+s_z=s\), the remaining graph is disconnected. Since transitivity satisfies the component maximum rule, we have
\[
F_x(s)=
\max_{\substack{s_y+s_z=s\\0\leq s_y\leq n_y\\0\leq s_z\leq n_z}}
\max\{F_y(s_y),F_z(s_z)\}.
\]
This is intentionally a recurrence for \(\Tr\), not a recurrence for \(D\).

Now suppose that \(x\) is a join node, so \(G_x=G_y\vee G_z\). The join case is the only place where parts may mix vertices from both sides. In a transitive partition of an induced subgraph of \(G_x\), call a part mixed if it contains vertices from both \(G_y\) and \(G_z\). Such a mixed part dominates every other part, so mixed parts can always be placed at the beginning of the ordering.

Let \(a\) be the number of vertices of \(G_y\) not used by the pure residual partition on the \(G_y\)-side, and let \(b\) be the analogous number for \(G_z\). These \(a+b\) vertices are reserved: some of them are truly deleted, and the others are used to form mixed parts. If the parent table entry deletes exactly \(s\) vertices in total, then \(L=a+b-s\) reserved vertices remain available for mixed parts.

Define \(q(s,a,b)\) as the maximum number of mixed parts that can be formed from these \(L\) remaining reserved vertices. If \(s>a+b\), the choice is infeasible. If \(L=0\), then \(q(s,a,b)=0\). If \(L=1\), the choice is infeasible because one vertex alone cannot form a mixed part. If \(L>0\) and either \(a=0\) or \(b=0\), the choice is also infeasible because a mixed part needs at least one vertex from each side. In the remaining case, where \(L\geq 2\) and \(a,b\geq 1\), we set \(q(s,a,b)=\min\{a,b,\lfloor L/2\rfloor\}.\) The term \(\min\{a,b\}\) says that every mixed part needs a vertex from each side, and the term \(\lfloor L/2\rfloor\) says that every mixed part has size at least two.

The join recurrence is
\[
F_x(s)=
\max_{\substack{0\leq a\leq n_y\\0\leq b\leq n_z\\(s,a,b)\text{ feasible}}}
\bigl(F_y(a)+F_z(b)+q(s,a,b)\bigr).
\]
The desired value is obtained at the root \(r\) of the cotree: \(D(G)=\Tr(G)=F_r(0)\).

\subsection{The dynamic programming algorithm}
\label{subsec:cograph-dp-algorithm}

We now present the dynamic programming algorithm. The binary cotree is processed in postorder, so when an internal node \(x\) is considered, the tables of its children are already available. For each node \(x\), we maintain a table \(F_x(0),F_x(1),\ldots,F_x(n_x)\), where \(F_x(s)\) denotes the maximum transitivity value obtainable from \(G_x\) after deleting exactly \(s\) vertices. Infeasible entries are assigned the value \(-\infty\).

If \(x\) is a leaf, then \(F_x(0)=1\) and \(F_x(1)=0\). If \(x\) is a union node with children \(y\) and \(z\), the algorithm tries all ways of distributing the deleted vertices between \(G_y\) and \(G_z\), and uses the component maximum rule for transitivity. If \(x\) is a join node, the algorithm guesses two integers \(a\) and \(b\), representing the reserved vertices on the two sides. These reserved vertices are not used in the pure transitive partitions of the children; they are either deleted or used to form mixed parts. A mixed part contains at least one vertex from each child, and hence it dominates all pure parts that appear later in the ordering. If \(s\) of the reserved vertices are deleted, then \(L=a+b-s\) vertices remain. The number of mixed parts is \(0\) when \(L=0\), and is \(\min\{a,b,\lfloor L/2\rfloor\}\) when \(L\geq2\) and \(a,b\geq1\); all other cases are infeasible. After processing all nodes, the desired value is \(F_r(0)\), where \(r\) is the root of the cotree, because no vertex is deleted from the original graph.

\begin{algorithm}[H]
	\caption{\textsc{Upper-Domatic-Cograph}\((G)\)}
	\label{alg:upper-domatic-cograph}
	\KwIn{A cograph \(G\)}
	\KwOut{The upper domatic number \(D(G)\)}
	
	Construct a binary cotree \(T\) of \(G\)\;
	Compute \(n_x=|V(G_x)|\) for every node \(x\) of \(T\)\;
	
	\ForEach{node \(x\) of \(T\), in postorder}{
		\eIf{\(x\) is a leaf}{
			\(F_x(0)\gets 1\), \(F_x(1)\gets 0\)\;
		}{
			Let \(y\) and \(z\) be the children of \(x\)\;
			Set \(F_x(s)\gets -\infty\) for all \(0\leq s\leq n_x\)\;
			
			\eIf{\(x\) is a union node}{
				\For{\(s_y\gets 0\) \KwTo \(n_y\)}{
					\For{\(s_z\gets 0\) \KwTo \(n_z\)}{
						\(s\gets s_y+s_z\)\;
						\(val\gets \max\{F_y(s_y),F_z(s_z)\}\)\;
						\(F_x(s)\gets \max\{F_x(s),val\}\)\;
					}
				}
			}{
				\tcp{\(x\) is a join node}
				\For{\(a\gets 0\) \KwTo \(n_y\)}{
					\For{\(b\gets 0\) \KwTo \(n_z\)}{
						\For{\(s\gets 0\) \KwTo \(a+b\)}{
							\(L\gets a+b-s\)\;
							
							\If{\(L=0\)}{
								\(q\gets 0\)\;
								\(val\gets F_y(a)+F_z(b)+q\)\;
								\(F_x(s)\gets \max\{F_x(s),val\}\)\;
							}
							
							\If{\(L\geq 2\) and \(a,b\geq 1\)}{
								\(q\gets \min\{a,b,\lfloor L/2\rfloor\}\)\;
								\(val\gets F_y(a)+F_z(b)+q\)\;
								\(F_x(s)\gets \max\{F_x(s),val\}\)\;
							}
						}
					}
				}
			}
		}
	}
	
	Let \(r\) be the root of \(T\)\;
	\Return{\(F_r(0)\)}\;
\end{algorithm}

\subsection{Correctness of the algorithm}
\label{subsec:cograph-correctness}

\begin{theorem}
	\label{thm:cograph-upper-domatic-correct}
	Algorithm~\ref{alg:upper-domatic-cograph} correctly computes the upper domatic number of a cograph \(G\).
\end{theorem}

\begin{proof}
	By Lemma~\ref{lem:cograph-D-equals-Tr}, it is enough to compute \(\Tr(G)\). We show that, for every cotree node \(x\) and every \(s\) with \(0\leq s\leq n_x\), the table entry \(F_x(s)\) computed by the algorithm is exactly the maximum value of \(\Tr(G_x-S)\) over all sets \(S\subseteq V(G_x)\) with \(|S|=s\).
	
	If \(x\) is a leaf, then \(G_x\) consists of one vertex. Hence, without deleting this vertex we have transitivity \(1\), and after deleting it we obtain the empty graph, whose transitivity is \(0\). Thus \(F_x(0)=1\) and \(F_x(1)=0\), as initialized by the algorithm.
	
	Now let \(x\) be an internal node with children \(y\) and \(z\). First, suppose that \(x\) is a union node. Then \(G_x=G_y\cup G_z\). If \(s_y\) vertices are deleted from \(G_y\) and \(s_z\) vertices are deleted from \(G_z\), with \(s_y+s_z=s\), the remaining graph is a disjoint union. Since the transitivity of a disconnected graph is the maximum transitivity of its components, the best value is obtained by maximizing \(\max\{F_y(s_y),F_z(s_z)\}\) over all such choices. This is exactly the update used by the algorithm at a union node.
	
	Now suppose that \(x\) is a join node, so \(G_x=G_y\vee G_z\). Consider an optimal transitive partition of \(G_x-S\), where \(|S|=s\). Let \(M_y\) and \(M_z\) be the sets of vertices from \(G_y\) and \(G_z\) that occur in mixed parts. Define \(A=(S\cap V(G_y))\cup M_y\) and \(B=(S\cap V(G_z))\cup M_z\), and let \(a=|A|\) and \(b=|B|\). The pure parts on the \(G_y\)-side form a transitive partition of \(G_y-A\), so they contribute at most \(F_y(a)\) parts. Similarly, the pure parts on the \(G_z\)-side contribute at most \(F_z(b)\) parts. The remaining \(a+b-s\) vertices are exactly the reserved vertices used in mixed parts. Since each mixed part needs at least one vertex from each side and at least two vertices in total, the number of mixed parts is at most \(q(s,a,b)\). Therefore every feasible transitive partition has size at most \(F_y(a)+F_z(b)+q(s,a,b)\) for some feasible triple \((s,a,b)\).
	
	Conversely, fix a feasible triple \((s,a,b)\). Choose sets \(A\subseteq V(G_y)\) and \(B\subseteq V(G_z)\) with \(|A|=a\) and \(|B|=b\) witnessing the values \(F_y(a)\) and \(F_z(b)\). If \(a+b=s\), then all vertices of \(A\cup B\) are deleted and no mixed part is formed. Otherwise, let \(L=a+b-s\). Feasibility gives \(L\geq 2\) and \(a,b\geq 1\). Let \(q=q(s,a,b)\). Since \(q\leq a\), \(q\leq b\), and \(2q\leq L\), we can choose \(q\) vertices from \(A\) and \(q\) vertices from \(B\), and pair them to form \(q\) mixed parts. We then retain \(L-2q\) further vertices from \(A\cup B\) and distribute them arbitrarily among these mixed parts. All remaining vertices of \(A\cup B\) are deleted. Thus, exactly \(s\) vertices are deleted and exactly \(q\) mixed parts are formed.
	
	Place these mixed parts first in the ordering. Then append a transitive partition of \(G_y-A\), followed by a transitive partition of \(G_z-B\). Every mixed part dominates every later part, because it contains vertices from both sides of the join. Also, every pure part on the \(G_y\)-side dominates every pure part on the \(G_z\)-side. Therefore, the resulting ordered partition is transitive and has \(F_y(a)+F_z(b)+q(s,a,b)\) parts. Hence, the join update is correct.
	
	By induction over the postorder traversal of the cotree, every table entry is computed correctly. At the root \(r\), no vertex is deleted from the original graph, so the algorithm returns \(F_r(0)=\Tr(G)\). By Lemma~\ref{lem:cograph-D-equals-Tr}, this is equal to \(D(G)\). Therefore, Algorithm~\ref{alg:upper-domatic-cograph} correctly computes the upper domatic number of \(G\).
\end{proof}

\subsection{Running time analysis}
\label{subsec:cograph-running-time}

\begin{theorem}
	\label{thm:cograph-upper-domatic}
	The upper domatic number of a cograph \(G\) on \(n\) vertices can be computed in \(O(n^3)\) time. If the input graph is not already given with a cotree, the cotree construction adds only \(O(n+m)\) time, and the total time remains \(O(n^3)\) for cographs.
\end{theorem}

\begin{proof}
	A leaf is processed in \(O(1)\) time. Let \(x\) be an internal node with children \(y\) and \(z\). If \(x\) is a union node, then the algorithm checks all pairs \((s_y,s_z)\), where \(0\leq s_y\leq n_y\) and \(0\leq s_z\leq n_z\). Hence the time spent at \(x\) is \(O((n_y+1)(n_z+1))\), which is \(O(n_y n_z+n_x)\).
	
	If \(x\) is a join node, then the algorithm checks all triples \((a,b,s)\), where \(0\leq a\leq n_y\), \(0\leq b\leq n_z\), and \(0\leq s\leq a+b\). For each fixed pair \((a,b)\), there are at most \(n_x+1\) possible values of \(s\), and each update takes constant time. Therefore, the time spent at a join node \(x\) is \(O(n_x n_y n_z)\).
	
	It remains to sum these bounds over the binary cotree. For every internal node \(x\), let \(y\) and \(z\) be its two children. The quantity \(n_y n_z\) counts the number of unordered pairs of leaves whose least common ancestor is \(x\). Hence, \(\displaystyle \sum_x n_y n_z=\binom{n}{2}\), where the sum is over all internal nodes of the binary cotree. Therefore, the total time spent at union nodes is at most \(O(n^2)\). For join nodes, since \(n_x\leq n\) for every node \(x\), we have \(\displaystyle \sum_x n_x n_y n_z\leq n\sum_x n_y n_z=O(n^3)\). The initialization of all tables contributes \(\displaystyle O\!\left(\sum_x n_x\right)=O(n^2)\). Thus, the dynamic program runs in \(O(n^3)\) time.
	
	A cotree of a cograph can be constructed in \(O(n+m)\) time~\cite{corneil1985linear}. Hence, the total running time is \(O(n^3+n+m)\), which is \(O(n^3)\) for cographs.
\end{proof}

\begin{remark}
	If an optimal partition is required, we store one back-pointer for each table entry. At a union node, the pointer records the pair \((s_y,s_z)\) and the child whose value is used. At a join node, the pointer records the triple \((a,b,s)\) and the number \(q(s,a,b)\). During reconstruction, the child solutions provide the two reserved sets \(A\) and \(B\). From \(a,b,s\), and \(q(s,a,b)\), we select the vertices that remain deleted and distribute the other reserved vertices among the mixed parts as in the constructive proof of the join recurrence. Thus, an optimal transitive partition, and hence an optimal upper domatic partition, can be reconstructed from the stored back-pointers.
\end{remark}


\section{Upper Domatic Number in Unicyclic Graphs}
\label{sec:unicyclic}

In this section, we study the upper domatic number of unicyclic graphs. Samuel and Joseph~\cite{samuel2020new} proved that \(D(G)=\Tr(G)\) for every unicyclic graph \(G\). Hence, for this graph class, computing the upper domatic number is equivalent to computing the transitivity. Our approach is based on a simple structural observation: deleting any edge of the unique cycle of a unicyclic graph produces a tree. Since the transitivity of a tree can be computed in linear time~\cite{hedetniemi2018transitivity,hedetniemi1982linear}, it remains only to decide whether adding back the deleted edge increases the transitivity by one.

For a nonnegative integer \(q\), write \([q]=\{1,2,\ldots,q\}\), where \([0]=\emptyset\). We use \(\mathbf{1}_{\mathcal P}\) for the indicator of a statement \(\mathcal P\). Let \(H\) be a graph. A map \(\lambda\colon V(H)\rightarrow [k]\) is called a \emph{transitive \(k\)-labelling} of \(H\) if \(\max\{\lambda(v):v\in V(H)\}=k\), and, for every vertex \(v\in V(H)\), the set \(\lambda(N_H(v))\) contains \([\lambda(v)-1]\). Thus, a vertex with label \(j\) has a neighbour with each label in \([j-1]\).

\begin{lemma}
	\label{lem:transitive-labelling-equivalence}
	A graph \(H\) has a transitive partition of order \(k\) if and only if it has a transitive \(k\)-labelling.
\end{lemma}

\begin{proof}
	Let \(\{V_1,V_2,\ldots,V_k\}\) be a transitive partition of \(H\), and assign label \(i\) to every vertex of \(V_i\). Since \(V_i\) dominates \(V_j\) whenever \(i<j\), every vertex labelled \(j\) has a neighbour with each label in \([j-1]\). Hence, the resulting map is a transitive \(k\)-labelling.
	
	Conversely, let \(\lambda\) be a transitive \(k\)-labelling, and put \(V_i=\{v\in V(H):\lambda(v)=i\}\). Since some vertex has label \(k\), and this vertex has a neighbour with every label in \([k-1]\), all the sets \(V_1,V_2,\ldots,V_k\) are nonempty. Moreover, if \(i<j\), then every vertex of \(V_j\) has a neighbour in \(V_i\). Therefore, \(V_i\) dominates \(V_j\), and \(\{V_1,V_2,\ldots,V_k\}\) is a transitive partition of \(H\).
\end{proof}

\subsection{Reduction to a Decision Problem on a Tree}
\label{subsec:unicyclic-one-edge-reduction}

We use the following edge-deletion and edge-addition result of Hedetniemi and Hedetniemi~\cite{hedetniemi2018transitivity}.

\begin{proposition}\cite{hedetniemi2018transitivity}
	\label{prop:one-edge-transitivity}
	Let \(H\) be a graph. If \(e\in E(H)\), then \(\Tr(H)-1\leq \Tr(H-e)\leq \Tr(H)\). If \(x\) and \(y\) are nonadjacent vertices of \(H\), then \(\Tr(H)\leq \Tr(H+xy)\leq \Tr(H)+1\).
\end{proposition}

Let \(G\) be a unicyclic graph, let \(xy\) be an arbitrary edge of its unique cycle, and let \(T=G-xy\). Then \(T\) is a tree. We compute \(\tau=\Tr(T)\) using the known linear-time tree algorithm. Since \(G=T+xy\), Proposition~\ref{prop:one-edge-transitivity} gives \(\tau\leq \Tr(G)\leq \tau+1\). Hence, \(\Tr(G)\in\{\tau,\tau+1\}\), and it remains only to decide whether \(G\) admits a transitive \((\tau+1)\)-labelling.

The following logarithmic bound ensures that the subset dynamic program described below runs in polynomial time.

\begin{lemma}
	\label{lem:tree-logarithmic-bound}
	If a tree \(T\) on \(n\) vertices has a transitive \(q\)-labelling, then \(n\geq 2^{q-1}\). Consequently, \(\Tr(T)\leq 1+\lfloor\log_2 n\rfloor\).
\end{lemma}

\begin{proof}
	For \(q\geq 1\), let \(m(q)\) be the minimum order of a tree having a transitive labelling in which some specified vertex has label \(q\). Clearly, \(m(1)=1\). Let \(v\) be a vertex labelled \(q\), where \(q\geq 2\). For each \(i\in[q-1]\), the vertex \(v\) has a neighbour \(v_i\) labelled \(i\). These neighbours are distinct, and they belong to distinct components of \(T-v\).
	
	Moreover, the restriction of the labelling to the component containing \(v_i\) still satisfies the transitive-labelling condition. Indeed, the only vertex of this component adjacent to \(v\) is \(v_i\), and the label \(q\) of \(v\) cannot supply any label in \([i-1]\). Therefore, the component containing \(v_i\) has at least \(m(i)\) vertices. It follows that \(m(q)\geq 1+\displaystyle\sum_{i=1}^{q-1}m(i)\). By induction, \(m(i)\geq 2^{i-1}\) for every \(i\), and hence \(m(q)\geq 1+\displaystyle\sum_{i=1}^{q-1}2^{i-1}=2^{q-1}\). The stated upper bound on \(\Tr(T)\) follows immediately.
\end{proof}

Set \(k=\tau+1\). By Lemma~\ref{lem:tree-logarithmic-bound}, we have \(k\leq 2+\lfloor\log_2 n\rfloor\). Thus, although the decision algorithm has an exponential dependence on \(k\), it is polynomial in \(n\).

\subsection{Testing Whether the Deleted Edge Increases Transitivity}
\label{subsec:unicyclic-edge-test}

Root the tree \(T\) at \(x\). For a vertex \(v\), let \(C(v)\) be its set of children, and let \(T_v\) be the subtree rooted at \(v\). We test all possible labels \(a,b\in[k]\) for the endpoints \(x\) and \(y\), respectively. For fixed \(a\) and \(b\), define
\[
\eta_{a,b}(v)=
\begin{cases}
	\{b\}, & \text{if \(v=x\)},\\
	\{a\}, & \text{if \(v=y\)},\\
	\emptyset, & \text{otherwise}.
\end{cases}
\]
The set \(\eta_{a,b}(v)\) records the label supplied to \(v\) by the deleted edge \(xy\) after that edge is restored.

For \(v\in V(T)\), \(c\in[k]\), \(p\in\{0,1,\ldots,k\}\), and \(h\in\{0,1\}\), we define a Boolean state \(F_v^{a,b}(c,p,h)\). The value \(p=0\) is used only for the root \(x\), and it means that the root has no parent. The state \(F_v^{a,b}(c,p,h)\) is true if and only if the vertices of \(T_v\) can be labelled from \([k]\) so that the following conditions hold:

\begin{enumerate}
	\item the vertex \(v\) has label \(c\);
	\item the fixed labels \(\lambda(x)=a\) and \(\lambda(y)=b\) are respected whenever the corresponding endpoint belongs to \(T_v\);
	\item every vertex of \(T_v\setminus\{v\}\) satisfies the transitive-labelling condition using its actual neighbours, together with any contribution of the restored edge specified by \(\eta_{a,b}\);
	\item the vertex \(v\) satisfies the transitive-labelling condition when it may additionally use an external parent of label \(p\); and
	\item \(h=1\) if and only if some vertex of \(T_v\) receives label \(k\).
\end{enumerate}

The recurrence is local. Fix a vertex \(v\) and a candidate label \(c\). For every child \(u\in C(v)\), choose a label \(d_u\in[k]\) and a value \(h_u\in\{0,1\}\) such that \(F_u^{a,b}(d_u,c,h_u)\) is true. Let
\(
M=\{d_u:u\in C(v),\,d_u<c\}
\)
and let
\(
h^*=\mathbf{1}_{\{c=k\}}\vee\bigvee_{u\in C(v)}h_u.
\)
Then \(F_v^{a,b}(c,p,h)=\mathrm{true}\) if and only if \(c=a\) whenever \(v=x\), \(c=b\) whenever \(v=y\), and there is a choice for the children such that \(h=h^*\) and
\[
[c-1]\subseteq M\cup\bigl(\{p\}\cap[k]\bigr)\cup\eta_{a,b}(v).
\tag{1}
\label{eq:unicyclic-dp-recurrence}
\]
Indeed, the labels in \(M\) are supplied by the children of \(v\), the label \(p\) is supplied by the parent of \(v\), and \(\eta_{a,b}(v)\) is supplied by the restored edge \(xy\).

The sets \(M\) are represented explicitly as subsets of \([c-1]\). For a fixed pair \((a,b)\), the states are computed in postorder. For each pair \((v,c)\), we begin with the reachable pair
\(
\left(\emptyset,\mathbf{1}_{\{c=k\}}\right).
\)
When a child \(u\) is processed, every current pair \((M,\widehat h)\) is combined with every pair \((d,h')\) such that \(F_u^{a,b}(d,c,h')=\mathrm{true}\), producing
\[
\left(M\cup(\{d\}\cap[c-1]),\,\widehat h\vee h'\right).\]

Each reachable family is maintained as a set, so duplicate pairs are discarded. For example, the pairs can be stored in a binary trie keyed by the characteristic vector of \(M\), followed by the Boolean coordinate. In this representation, inserting a pair and checking whether it is already present takes \(O(c)\) time. After all children have been processed, Equation~\eqref{eq:unicyclic-dp-recurrence} determines the states \(F_v^{a,b}(c,p,h)\). The pair \((a,b)\) is feasible precisely when \(F_x^{a,b}(a,0,1)=\mathrm{true}\).

\subsection{The decision subroutine and the complete algorithm}
\label{subsec:unicyclic-algorithm}

Before presenting the main algorithm, we describe the decision subroutine that tests whether adding back the deleted cycle edge increases the transitivity. Recall that \(T=G-xy\) is a tree, and we want to decide whether \(T+xy\) admits a transitive partition of order \(k\). The subroutine tries all possible labels \(a,b\in[k]\) for the endpoints \(x\) and \(y\) of the deleted edge. For each fixed pair \((a,b)\), it roots \(T\) at \(x\) and processes the vertices in postorder, so that the information for every child subtree is already available before its parent is processed. The state \(F_v^{a,b}(c,p,h)\) records whether the subtree \(T_v\) admits a consistent transitive \(k\)-labelling under the assumptions that \(v\) receives label \(c\), the parent of \(v\) has label \(p\), and the Boolean value \(h\) indicates whether label \(k\) appears somewhere in \(T_v\). The only edge of \(G\) not present in \(T\) is \(xy\), and its possible contribution is handled separately through the set \(\eta_{a,b}(v)\). Thus, the dynamic program works entirely on the rooted tree \(T\), while still accounting for the effect of restoring the edge \(xy\).

\begin{algorithm}[H]
	\small
	\caption{\textsc{One-Edge-Increase-Test}\((T,x,y,k)\)}
	\label{alg:one-edge-increase-test}
	\KwIn{A tree \(T\), two vertices \(x,y\in V(T)\), and an integer \(k\)}
	\KwOut{\(\mathrm{true}\) if \(T+xy\) has a transitive partition of order \(k\); otherwise \(\mathrm{false}\)}
	
	Root \(T\) at \(x\) and compute a postorder of its vertices\;
	
	\For{\(a\gets 1\) \KwTo \(k\)}{
		\For{\(b\gets 1\) \KwTo \(k\)}{
			Initialize every state \(F_v^{a,b}(c,p,h)\) to \(\mathrm{false}\)\;
			
			\ForEach{vertex \(v\) in postorder}{
				\For{\(c\gets 1\) \KwTo \(k\)}{
					\If{\((v=x\text{ and }c\neq a)\) or \((v=y\text{ and }c\neq b)\)}{
						continue\;
					}
					
					Set \(\mathcal R\gets\{(\emptyset,\mathbf{1}_{\{c=k\}})\}\)\;
					
					\ForEach{\(u\in C(v)\)}{
						Set \(\mathcal R'\gets\emptyset\)\;
						
						\ForEach{\((M,\widehat h)\in\mathcal R\)}{
							\ForEach{\(d\in[k]\) and \(h'\in\{0,1\}\) with \(F_u^{a,b}(d,c,h')=\mathrm{true}\)}{
								Insert \(\left(M\cup(\{d\}\cap[c-1]),\widehat h\vee h'\right)\) into \(\mathcal R'\) if it is not already present\;
							}
						}
						
						Set \(\mathcal R\gets\mathcal R'\)\;
					}
					
					Let \(P_v=\{0\}\) if \(v=x\), and let \(P_v=[k]\) otherwise\;
					
					\ForEach{\(p\in P_v\) and \(g\in\{0,1\}\)}{
						\If{there exists \((M,g)\in\mathcal R\) satisfying Equation~\eqref{eq:unicyclic-dp-recurrence}}{
							Set \(F_v^{a,b}(c,p,g)\gets\mathrm{true}\)\;
						}
					}
				}
			}
			
			\If{\(F_x^{a,b}(a,0,1)=\mathrm{true}\)}{
				\Return{\(\mathrm{true}\)}\;
			}
		}
	}
	
	\Return{\(\mathrm{false}\)}\;
\end{algorithm}

The complete algorithm now follows directly. We delete one edge \(xy\) of the unique cycle, compute the transitivity \(\tau\) of the resulting tree \(T=G-xy\), and then use the above subroutine to decide whether \(G=T+xy\) has a transitive partition of order \(k=\tau+1\).

\begin{algorithm}[H]
	\small
	\caption{\textsc{Transitivity-Unicyclic}\((G)\)}
	\label{alg:transitivity-unicyclic-new}
	\KwIn{A unicyclic graph \(G\)}
	\KwOut{The transitivity \(\Tr(G)\)}
	
	Find the unique cycle of \(G\) and choose an arbitrary cycle edge \(xy\)\;
	Set \(T\gets G-xy\)\;
	Compute \(\tau\gets\Tr(T)\) using the tree algorithm\;
	Set \(k\gets\tau+1\)\;
	
	\eIf{\textsc{One-Edge-Increase-Test}\((T,x,y,k)\) returns \(\mathrm{true}\)}{
		\Return{\(k\)}\;
	}{
		\Return{\(\tau\)}\;
	}
\end{algorithm}

\subsection{Correctness}
\label{subsec:unicyclic-new-correctness}

\begin{lemma}
	\label{lem:unicyclic-dp-state-correctness}
	For fixed endpoint labels \(a,b\in[k]\), the state \(F_v^{a,b}(c,p,h)\) is \(\mathrm{true}\) if and only if a labelling satisfying the conditions in its definition exists for the subtree \(T_v\).
\end{lemma}

\begin{proof}
	We prove the lemma by induction on the number of vertices of \(T_v\). If \(v\) is a leaf, then it has no child. Hence, the only possible contributions to the labels required by \(v\) come from its parent and, when \(v\in\{x,y\}\), from the restored edge \(xy\). Thus, Equation~\eqref{eq:unicyclic-dp-recurrence} is exactly the transitive-labelling condition at \(v\), and the statement follows.
	
	Now suppose that \(v\) has at least one child and that the statement holds for every child subtree. First assume that the recurrence declares \(F_v^{a,b}(c,p,h)\) to be \(\mathrm{true}\). For every child \(u\in C(v)\), choose a pair \((d_u,h_u)\) used by the recurrence. By the induction hypothesis, \(T_u\) has a labelling represented by \(F_u^{a,b}(d_u,c,h_u)\). The child subtrees are pairwise vertex-disjoint, and the only tree edge joining \(T_u\) to the rest of \(T_v\) is \(uv\). The restored edge \(xy\), whenever it involves a vertex of \(T_u\), is already accounted for by the contribution \(\eta_{a,b}\). Hence, the labellings of the child subtrees can be combined, and \(v\) can be assigned label \(c\).
	
	Every vertex in a child subtree satisfies the required condition by the induction hypothesis. For the vertex \(v\), Equation~\eqref{eq:unicyclic-dp-recurrence} guarantees that every label in \([c-1]\) appears on a neighbour of \(v\), supplied either by a child, by the parent of \(v\), or by the restored edge \(xy\). Thus, the combined labelling satisfies all conditions in the definition of \(F_v^{a,b}(c,p,h)\), and the value of \(h\) is correct.
	
	Conversely, suppose that a labelling satisfying the conditions in the definition of \(F_v^{a,b}(c,p,h)\) exists for \(T_v\). Let \(d_u\) be the label of each child \(u\), and let \(h_u\) indicate whether label \(k\) occurs in \(T_u\). By the induction hypothesis, \(F_u^{a,b}(d_u,c,h_u)\) is \(\mathrm{true}\). Every lower label required by \(v\) is supplied either by one of its children, by its parent, or by the restored edge \(xy\). Therefore, the corresponding set \(M\) satisfies Equation~\eqref{eq:unicyclic-dp-recurrence}, and the recurrence declares \(F_v^{a,b}(c,p,h)\) to be \(\mathrm{true}\). This completes the induction.
\end{proof}

\begin{lemma}
	\label{lem:one-edge-test-correctness}
	Algorithm~\ref{alg:one-edge-increase-test} returns \(\mathrm{true}\) if and only if \(T+xy\) has a transitive partition of order \(k\).
\end{lemma}

\begin{proof}
	Suppose that the algorithm returns \(\mathrm{true}\). Then, for some \(a,b\in[k]\), the state \(F_x^{a,b}(a,0,1)\) is \(\mathrm{true}\). By Lemma~\ref{lem:unicyclic-dp-state-correctness}, the whole tree \(T=T_x\) has a labelling satisfying the transitive-labelling condition after the edge \(xy\) is restored. Since the final state has value \(h=1\), some vertex receives label \(k\). Hence, the labelling is a transitive \(k\)-labelling of \(T+xy\), and Lemma~\ref{lem:transitive-labelling-equivalence} gives a transitive partition of order \(k\).
	
	Conversely, suppose that \(T+xy\) has a transitive partition of order \(k\), and let \(\lambda\) be the corresponding transitive \(k\)-labelling. Set \(a=\lambda(x)\) and \(b=\lambda(y)\). Applying Lemma~\ref{lem:unicyclic-dp-state-correctness} bottom-up to the restrictions of \(\lambda\) to the rooted subtrees gives \(F_x^{a,b}(a,0,1)=\mathrm{true}\). Thus, the algorithm returns \(\mathrm{true}\) for the pair \((a,b)\).
\end{proof}

\begin{theorem}
	\label{thm:unicyclic-new-correctness}
	Algorithm~\ref{alg:transitivity-unicyclic-new} correctly computes the transitivity of a unicyclic graph.
\end{theorem}

\begin{proof}
	Let \(xy\) be the chosen cycle edge, let \(T=G-xy\), and let \(\tau=\Tr(T)\). By Proposition~\ref{prop:one-edge-transitivity}, we have \(\Tr(G)\in\{\tau,\tau+1\}\). By Lemma~\ref{lem:one-edge-test-correctness}, Algorithm~\ref{alg:one-edge-increase-test} correctly determines whether \(G=T+xy\) has a transitive partition of order \(\tau+1\). Therefore, Algorithm~\ref{alg:transitivity-unicyclic-new} returns \(\tau+1\) exactly when \(\Tr(G)=\tau+1\), and otherwise it returns \(\tau\). Hence, it always returns \(\Tr(G)\).
\end{proof}

Since Samuel and Joseph~\cite{samuel2020new} proved that \(D(G)=\Tr(G)\) for every unicyclic graph, we obtain the following consequence.

\begin{corollary}
	\label{cor:unicyclic-new-upper-domatic}
	Algorithm~\ref{alg:transitivity-unicyclic-new} correctly computes the upper domatic number of a unicyclic graph.
\end{corollary}

\subsection{Running Time}
\label{subsec:unicyclic-new-running-time}

\begin{theorem}
	\label{thm:unicyclic-new-running-time}
	Let \(G\) be a unicyclic graph on \(n\) vertices. Algorithm~\ref{alg:transitivity-unicyclic-new} runs in \(O(n^3)\) time.
\end{theorem}

\begin{proof}
	The unique cycle of \(G\) can be found in \(O(n)\) time, and after deleting the chosen cycle edge \(xy\), the graph \(T=G-xy\) is a tree. The value \(\tau=\Tr(T)\) can be computed in \(O(n)\) time using the known tree algorithm. Set \(k=\tau+1\).
	
	We now estimate the running time of Algorithm~\ref{alg:one-edge-increase-test}. The algorithm tests all \(k^2\) possible choices for the labels \(a\) and \(b\) of the endpoints \(x\) and \(y\). Fix one such pair \((a,b)\).
	
	For a fixed vertex \(v\) and a fixed label \(c\), the reachable family \(\mathcal R\) consists of pairs \((M,h)\), where \(M\subseteq[c-1]\) and \(h\in\{0,1\}\). Therefore, \(\mathcal R\) contains at most \(2^c\) distinct pairs. Each pair is stored as a \(0\)-\(1\) string of length \(c\), namely the characteristic vector of \(M\) followed by the bit \(h\). We keep these strings in a binary trie, so inserting a pair and checking whether it is already present takes \(O(c)\) time. For a fixed child \(u\), every pair in \(\mathcal R\) is combined with at most \(2k\) choices \((d,h')\), where \(d\in[k]\) and \(h'\in\{0,1\}\). Hence, for a fixed pair \((a,b)\), all child transitions take \(\displaystyle \sum_{vu\in E(T)}\sum_{c=1}^{k}O(kc2^c)=O(nk^2 2^k)\) time.
	
	After all children of \(v\) have been processed, the algorithm checks the possible parent labels \(p\) and Boolean values \(g\). For fixed \(v\) and \(c\), there are at most \(2k\) such choices. Scanning the reachable family and checking Equation~\eqref{eq:unicyclic-dp-recurrence} takes \(O(c2^c)\) time for each such choice. Hence, over all vertices and all values of \(c\), these final checks also take \(O(nk^2 2^k)\) time. Thus, the total running time for a fixed pair \((a,b)\) is \(O(nk^2 2^k)\).
	
	Since there are \(k^2\) choices for \((a,b)\), Algorithm~\ref{alg:one-edge-increase-test} runs in \(O(nk^4 2^k)\) time. By Lemma~\ref{lem:tree-logarithmic-bound},
	\(
	k\leq 2+\lfloor\log_2 n\rfloor,
	\)
	and hence \(2^k\leq 4n\). Moreover, \(k=O(\log n)\), and \((\log n)^4=O(n)\); therefore, \(k^4=O(n)\). It follows that
	\(
	O(nk^4 2^k)
	=O(n\cdot n\cdot n)
	=O(n^3).
	\)
	The initial \(O(n)\)-time steps are dominated by this bound. Therefore, Algorithm~\ref{alg:transitivity-unicyclic-new} runs in \(O(n^3)\) time.
\end{proof}

\begin{theorem}
	\label{thm:unicyclic-final}
	The transitivity and the upper domatic number of a unicyclic graph on \(n\) vertices can be computed in \(O(n^3)\) time.
\end{theorem}


\section{Upper Domatic Number in the Complement of Bipartite Graphs}
\label{sec:complement-bipartite}

In this section, we prove that the upper domatic number and the transitivity are equal for complements of bipartite graphs. We first show that every such graph admits a maximum upper domatic partition containing a source set.

\begin{lemma}
	\label{upper_domatic_partition_CBG}
	Let \(G\) be the complement of a bipartite graph with vertex partition \(V(G)=X\cup Y\), where \(X\) and \(Y\) are cliques in \(G\). If \(D(G)=k\), then \(G\) has an upper domatic partition \(\pi=\{V_1,V_2,\ldots,V_k\}\) containing a source set.
\end{lemma}

\begin{proof}
	Let \(\pi=\{V_1,V_2,\ldots,V_k\}\) be a \(D\)-partition of \(G\). Since \(G\) is the complement of a bipartite graph, both \(X\) and \(Y\) are cliques in \(G\). If all parts of \(\pi\) are singletons, then the upper domatic condition implies that every two vertices of \(G\) are adjacent. Hence \(G\) is complete, and any singleton part is a source set. Therefore, assume that some part \(V_i\) has size at least two.
	
	If \(V_i\) contains vertices from both \(X\) and \(Y\), then \(V_i\) dominates every other part. Indeed, a vertex of \(V_i\cap X\) dominates all vertices of \(X\setminus V_i\), and a vertex of \(V_i\cap Y\) dominates all vertices of \(Y\setminus V_i\). Thus \(V_i\) is a source set. Hence, we may assume that every part of \(\pi\) is contained entirely in \(X\) or entirely in \(Y\).
	
	Without loss of generality, let \(V_i=\{x_{t_1},x_{t_2},\ldots,x_{t_s}\}\subseteq X\), where \(s\geq 2\). If \(Y=\emptyset\), then \(V_i\) dominates every other part because \(G[X]\) is complete, and hence \(V_i\) is a source set. Now suppose that \(Y\neq\emptyset\). The part \(V_i\) dominates every part contained in \(X\). Moreover, let \(V_j=\{y\}\) be a singleton part contained in \(Y\). Since \(\pi\) is an upper domatic partition, either \(V_i\to V_j\) or \(V_j\to V_i\). In either case, \(y\) has a neighbour in \(V_i\). Hence \(V_i\to V_j\). Thus \(V_i\) dominates every singleton part contained in \(Y\). Therefore, if every part contained in \(Y\) is a singleton, then \(V_i\) is a source set.
	
	It remains to consider the case where some part contained in \(Y\) has size at least two. Choose such a part \(V_j=\{y_{p_1},y_{p_2},\ldots,y_{p_r}\}\subseteq Y\), where \(r\geq 2\). If \(V_i\) dominates every part contained in \(Y\) of size at least two, then, together with the previous paragraph, \(V_i\) dominates every other part of \(\pi\), and hence \(V_i\) is a source set.
	
	Therefore, assume that \(V_i\) does not dominate some such part \(V_j\). Since \(\pi\) is an upper domatic partition, \(V_j\) must dominate \(V_i\). We construct a new partition \(\pi'\) by exchanging \(x_{t_1}\) and \(y_{p_1}\): replace \(x_{t_1}\) in \(V_i\) by \(y_{p_1}\), replace \(y_{p_1}\) in \(V_j\) by \(x_{t_1}\), and keep every other part unchanged. Since \(s\geq 2\) and \(r\geq 2\), both modified parts are mixed; that is, each contains at least one vertex from \(X\) and at least one vertex from \(Y\). Therefore, each of these two modified parts dominates every other part of \(\pi'\). All domination relations between unchanged parts remain the same. Hence \(\pi'\) is an upper domatic partition of size \(k\), and the modified part corresponding to \(V_i\) is a source set.
\end{proof}

The previous lemma shows that complements of bipartite graphs always admit a maximum upper domatic partition with a source set. This allows us to order one part first and then apply induction to the remaining graph. Using this idea, we prove that, in this class, the upper domatic number and transitivity are equal.

\begin{theorem}
	\label{UDN_equal_tranitivity_theorem_BCG}
	Let \(G\) be the complement of a bipartite graph. Then \(\Tr(G)=D(G)\).
\end{theorem}

\begin{proof}
	For every graph \(G\), every transitive partition is also an upper domatic partition. Hence \(\Tr(G)\leq D(G)\). It remains to prove the reverse inequality for complements of bipartite graphs.
	
	Suppose, to the contrary, that there exists a complement of a bipartite graph \(G_0\) with \(\Tr(G_0)<D(G_0)\). Choose such a graph with the minimum number of vertices. Let \(D(G_0)=k\), and let \(\pi=\{V_1,V_2,\ldots,V_k\}\) be a \(D\)-partition of \(G_0\). By Lemma~\ref{upper_domatic_partition_CBG}, we may assume that \(V_1\) is a source set. Since \(k=1\) would imply \(\Tr(G_0)=D(G_0)=1\), we have \(k\geq 2\).
	
	Let \(G'=G_0-V_1\), and let \(G_1,G_2,\ldots,G_s\) be the connected components of \(G'\). Each \(G_i\) is again the complement of a bipartite graph and has fewer vertices than \(G_0\). Hence, by the minimality of \(G_0\), we have \(\Tr(G_i)=D(G_i)\) for every \(1\leq i\leq s\).
	
	By the component property of the upper domatic number, we have \(D(G')\leq \max\{D(G_i):1\leq i\leq s\}\) by~\cite{samuel2020new}. Also, by the component property of transitivity, we have \(\Tr(G')=\max\{\Tr(G_i):1\leq i\leq s\}\) by~\cite{hedetniemi2018transitivity}. Since \(\Tr(G_i)=D(G_i)\) for every \(i\), it follows that \(\Tr(G')\geq D(G')\).
	
	Since \(V_1\) is a source set in \(\pi\), deleting \(V_1\) from \(\pi\) gives an upper domatic partition of \(G'\) of size \(D(G_0)-1\). Thus \(D(G')\geq D(G_0)-1\). On the other hand, \(V_1\) dominates every vertex of \(G'\). Hence, placing \(V_1\) before a maximum transitive partition of \(G'\) gives a transitive partition of \(G_0\). Therefore, \(\Tr(G_0)\geq \Tr(G')+1\geq D(G')+1\geq D(G_0)\), which contradicts \(\Tr(G_0)<D(G_0)\). Hence \(\Tr(G)=D(G)\) for every complement of a bipartite graph \(G\).
\end{proof}

It was proved in~\cite{paulSantra2025transitivityGraphClasses} that the transitivity of complements of bipartite chain graphs can be computed in linear time. By Theorem~\ref{UDN_equal_tranitivity_theorem_BCG}, the same holds for the upper domatic number.

\begin{corollary}
	\label{cor:upper-domatic-complement-chain}
	The upper domatic number of the complement of a bipartite chain graph can be computed in linear time.
\end{corollary}

\section{Upper Domatic Number in Split Graphs}
\label{sec:split-graphs}

In this section, we show that the upper domatic number of a split graph can be computed in linear time. Throughout this section, let \(G=(S\cup K,E)\) be a split graph, where \(S\) is an independent set and \(K\) is a clique of \(G\). We assume that \(K\) is a maximum clique of \(G\), that is, \(\omega(G)=|K|\). We first prove that the upper domatic number of \(G\) is either \(\omega(G)\) or \(\omega(G)+1\).

\begin{lemma}
	\label{split_graph_upper_domatic_lemma_1}
	Let \(G=(S\cup K,E)\) be a split graph with \(K\) a maximum clique. Then \(\omega(G)\leq D(G)\leq \omega(G)+1\).
\end{lemma}

\begin{proof}
	Let \(K=\{x_1,x_2,\ldots,x_{\omega(G)}\}\). Now, construct a vertex partition \(\pi=\{V_1, V_2, \ldots,V_{\omega(G)}\}\) as follows. Put \(V_1=S\cup\{x_1\}\), and for every \(2\leq i\leq \omega(G)\), put \(V_i=\{x_i\}\). Since \(K\) is a clique, \(V_i\to V_j\) for every \(1\leq i<j\leq \omega(G)\). Hence \(\pi\) is a transitive partition, and therefore also an upper domatic partition, of size \(\omega(G)\). Thus \(D(G)\geq \omega(G)\).
	
	Now we prove that \(D(G)\leq \omega(G)+1\). Suppose, to the contrary, that \(D(G)\geq \omega(G)+2\), and let \(\pi=\{V_1,V_2,\ldots,V_{\omega(G)+2}\}\) be an upper domatic partition of \(G\) of size \(\omega(G)+2\). Since \(|K|=\omega(G)\), at least two parts of \(\pi\), say \(V_i\) and \(V_j\), contain no vertices from \(K\). Hence both \(V_i\) and \(V_j\) are contained in \(S\). Since \(S\) is an independent set, neither \(V_i\) dominates \(V_j\) nor \(V_j\) dominates \(V_i\). This contradicts the fact that \(\pi\) is an upper domatic partition. Therefore, \(\omega(G)\leq D(G)\leq \omega(G)+1\).
\end{proof}

\begin{theorem}
	\label{split_graph_upper_domatic_lemma_2}
	Let \(G=(S\cup K,E)\) be a split graph with \(K\) a maximum clique. Then \(D(G)=\Tr(G)\). More precisely, if every vertex of \(K\) has a neighbour in \(S\), then \(D(G)=\Tr(G)=\omega(G)+1\); otherwise, \(D(G)=\Tr(G)=\omega(G)\).
\end{theorem}

\begin{proof}
	We divide the proof into two cases.
	
	First, suppose that there exists a vertex \(x\in K\) such that \(x\) has no neighbour in \(S\). We show that \(D(G)\neq \omega(G)+1\). Suppose, to the contrary, that \(D(G)=\omega(G)+1\). Let \(\pi=\{V_1,V_2,\ldots,V_{\omega(G)+1}\}\) be an upper domatic partition of \(G\) of size \(\omega(G)+1\). Since \(K\) has exactly \(\omega(G)\) vertices, at least one part of \(\pi\) contains no vertex from \(K\). Since two parts contained in \(S\) cannot be comparable under domination, exactly one part of \(\pi\) is contained in \(S\), and every other part contains exactly one vertex from \(K\).
	
	Let \(V_i\) be the part contained in \(S\), and suppose \(x\in V_j\), where \(j\neq i\). Since \(x\) has no neighbour in \(S\), and since \(S\) is independent, there is no edge between \(V_i\) and \(V_j\). Hence neither \(V_i\) dominates \(V_j\) nor \(V_j\) dominates \(V_i\), contradicting the assumption that \(\pi\) is an upper domatic partition. Therefore, \(D(G)\neq \omega(G)+1\). By Lemma~\ref{split_graph_upper_domatic_lemma_1}, we have \(D(G)=\omega(G)\). The construction in Lemma~\ref{split_graph_upper_domatic_lemma_1} gives a transitive partition of size \(\omega(G)\). Hence \(\omega(G)\leq \Tr(G)\leq D(G)=\omega(G)\), and so \(D(G)=\Tr(G)=\omega(G)\).
	
	Now suppose that every vertex of \(K\) has a neighbour in \(S\). Define a partition \(\pi=\{V_1,V_2,\ldots,V_{\omega(G)+1}\}\) by setting \(V_1=S\), and by assigning each vertex of \(K\) to a distinct singleton part among \(V_2,\ldots,V_{\omega(G)+1}\). Since every vertex of \(K\) has a neighbour in \(S\), the part \(V_1\) dominates every singleton part in \(K\). Since \(K\) is a clique, the singleton parts in \(K\) dominate each other in the prescribed order. Thus \(\pi\) is a transitive partition of size \(\omega(G)+1\), and hence \(\Tr(G)\geq \omega(G)+1\). By Lemma~\ref{split_graph_upper_domatic_lemma_1}, \(D(G)\leq \omega(G)+1\). Since \(\Tr(G)\leq D(G)\) holds for every graph, we get \(D(G)=\Tr(G)=\omega(G)+1\).
	
	Therefore, in both cases, \(D(G)=\Tr(G)\), with the claimed value.
\end{proof}

In~\cite{paulSantra2025transitivityGraphClasses}, it was proved that the transitivity of split graphs can be computed in linear time. Moreover, by Theorem~\ref{split_graph_upper_domatic_lemma_2}, \(D(G)=\Tr(G)\) for split graphs. Therefore, we obtain the following corollary.

\begin{corollary}
	\label{cor:upper-domatic-split}
	The upper domatic number of a split graph can be computed in linear time.
\end{corollary}

\section{Conclusion}
\label{sec:conclusion}

In this paper, we studied the upper domatic number from the viewpoints of computational hardness and exact algorithms on graph classes. We first proved that the decision problem is NP-complete by a reduction from \textsc{Clique}. The construction is simple but robust: it produces graphs obtained by joining an arbitrary graph to an independent false-twin module. As a consequence, the problem remains NP-complete even on connected graphs of diameter two. This answers the general variable-threshold version of the complexity question raised in Phillips's thesis.

On the algorithmic side, we identified several classes where the upper domatic number can be computed efficiently. A recurring theme is the relationship between the upper domatic number and transitivity. For cographs, we proved \(D(G)=\Tr(G)\), despite the fact that \(D(G)\) does not satisfy a simple equality rule over disconnected components. This allowed us to compute \(D(G)\) by a cotree dynamic program in \(O(n^3)\) time. For unicyclic graphs, we used the known equality \(D(G)=\Tr(G)\) and designed an \(O(n^3)\)-time algorithm for computing transitivity. The algorithm deletes an edge of the unique cycle, computes the transitivity of the resulting tree, and then uses a dynamic program to decide whether adding back the deleted edge increases the transitivity by one. For complements of bipartite graphs, we proved that every graph in the class admits a maximum upper domatic partition with a source set, which implies \(D(G)=\Tr(G)\). This gives a linear-time algorithm for complements of bipartite chain graphs. Finally, for split graphs, we showed that \(D(G)=\Tr(G)\) and that \(D(G)\) is either \(\omega(G)\) or \(\omega(G)+1\), yielding a linear-time algorithm.

Several natural questions remain open. The NP-completeness result obtained here applies even to connected diameter-two graphs, but it does not settle the complexity on more structured classes such as chordal graphs, bipartite graphs, or block graphs. It would also be interesting to improve the running time of the algorithm for unicyclic graphs further, or to find a more direct formula for the upper domatic number on that class. More generally, the equality \(D(G)=\Tr(G)\) plays a central role in the positive results of this paper. Characterizing graph classes for which this equality holds, or identifying structural obstructions to it, appears to be a promising direction for future work.

\bibliographystyle{plain}
\bibliography{UDN_Bib}

\end{document}